\documentclass[10 pt]{amsart}

\usepackage{etoolbox}
\makeatletter
\let\ams@starttoc\@starttoc
\makeatother
\usepackage[parfill]{parskip}
\makeatletter
\let\@starttoc\ams@starttoc
\patchcmd{\@starttoc}{\makeatletter}{\makeatletter\parskip\z@}{}{}
\makeatother

\usepackage[colorlinks=true,linkcolor=blue,citecolor=blue,urlcolor=blue]{hyperref}
\usepackage{bookmark}
\usepackage{amsthm,thmtools,amssymb,amsmath}
\usepackage{esint}
\usepackage{enumerate}

\usepackage{pictexwd,dcpic}

\swapnumbers

\numberwithin{equation}{section}

\usepackage{cleveref}
\crefname{lemma}{Lemma}{Lemmata}
\crefname{prop}{Proposition}{Propositions}
\crefname{thm}{Theorem}{Theorems}
\crefname{cor}{Corollary}{Corollaries}
\crefname{defn}{Definition}{Definitions}
\crefname{example}{Example}{Examples}
\crefname{rem}{Remark}{Remarks}
\crefname{ass}{Assumption}{Assumptions}
\crefname{not}{Notation}{Notation}

\renewcommand{\(}{\left(}
\renewcommand{\)}{\right)}
\renewcommand{\~}{\tilde}
\renewcommand{\-}{\bar}

\newcommand{\cn}{\colon}
\newcommand{\sub}{\subset}

\newcommand{\N}{\mathbb{N}}
\newcommand{\R}{\mathbb{R}}

\renewcommand{\S}{\mathbb{S}}
\renewcommand{\H}{\mathbb{H}}

\renewcommand{\a}{\alpha}
\renewcommand{\b}{\beta}

\renewcommand{\d}{\delta}
\newcommand{\e}{\epsilon}
\renewcommand{\k}{\kappa}
\renewcommand{\l}{\lambda}

\newcommand{\s}{\sigma}

\newcommand{\vt}{\vartheta}
\renewcommand{\O}{\Omega}

\newcommand{\inpr}[2]{\langle #1,#2 \rangle}
\newcommand{\fr}[2]{\frac{#1}{#2}}

\DeclareMathOperator{\id}{id}

\DeclareMathOperator{\osc}{osc}

\DeclareMathOperator{\dist}{dist}

\newcommand{\Theo}[3]{\begin{#1}\label{#2} #3 \end{#1}}
\newcommand{\pf}[1]{\begin{proof} #1 \end{proof}}
\newcommand{\eq}[1]{\begin{equation}\begin{alignedat}{2} #1 \end{alignedat}\end{equation}}

\newcommand{\ra}{\rightarrow}
\newcommand{\hra}{\hookrightarrow}

\newcommand{\mc}{\mathcal}
\renewcommand{\it}{\textit}
\newcommand{\mrm}{\mathrm}

\protected\def\ignorethis#1\endignorethis{}
\let\endignorethis\relax

\DeclareMathOperator{\argmin}{argmin}
\begin{document}

\title{Explicit rigidity of almost-umbilical hypersurfaces}


\author{Julien Roth and Julian Scheuer}
\date{\today}

\begin{abstract}
We give an explicit estimate of the distance of a closed, connected, oriented and immersed hypersurface of a space form to a geodesic sphere and show that the spherical closeness can be controlled by a power of an integral norm of the traceless second fundamental form, whenever the latter is sufficiently small. Furthermore we use the inverse mean curvature flow in the hyperbolic space to deduce the best possible order of decay in the class of $C^{\infty}$-bounded hypersurfaces of the Euclidean space.\end{abstract}

\address{Julien Roth, Laboratoire d'Analyse et de Math\'ematiques Appliqu\'ees, UPEM-UPEC, CNRS, F-77454 Marne-la-Vall\'ee, France}
\email{julien.roth@u-pem.fr}

\address{Julian Scheuer, Albert-Ludwigs-Universit{\"a}t, Mathematisches Institut, Eckerstr.~1,
79104 Freiburg, Germany}
\email{julian.scheuer@math.uni-freiburg.de}

\maketitle

\tableofcontents

\section{Introduction}

In this paper we prove two stability theorems of almost-umbilicity type, which give an answer to a question raised in \cite{Perez:/2011} and thereby partially improve \cite[Thm.~1.3, Thm.~1.4]{GrosjeanRoth:/2012}. Furthermore we use a recent counterexample for the inverse mean curvature flow in the hyperbolic space, cf.~\cite{HungWang:12/2014}, to provide a new counterexample for spherical closeness estimates.

 Let us shortly introduce the relevant notation. For an oriented hypersurface of a Riemannian manifold, $M^{n}\hra N^{n+1},$ $g$ denotes its induced metric, $|M|$ its surface area, $A$ its second fundamental form, $\mathring{A}$ the traceless part of $A,$
\eq{\mathring{A}=A-Hg,}
$x_{M}$ the center of mass of $M$ and $d_{\mc{H}}$ the Hausdorff distance of sets.

For a tensor field $(T_{i_{1}\dots i_{k}}^{j_{1}\dots j_{l}})$ on $M,$ we define its $L^{p}$-norm to be
\eq{\label{Lp}\|T\|_{p}=\(\int_{M}|T_{i_{1}\dots i_{k}}^{j_{1}\dots j_{l}}T^{i_{1}\dots i_{k}}_{j_{1}\dots j_{l}}|^{\fr p2}\)^{\fr 1p},}
where indices are raised or lowered with the help of $g.$ Let us formulate our first main result.

\Theo{thm}{Main}{
Let $M\hra\R^{n+1}$ be a closed, connected, oriented and immersed $C^{2}$-hypersurface with $|M|=1.$ Let $p>n\geq 2.$ Then there exist constants $c,\e_{0}>0$ depending on $n, p$ and $\|A\|_{p},$ as well as a constant $\a=\a(n,p),$ such that whenever there holds
\eq{\label{Main2}\|\mathring{A}\|_{p}< \|H\|_{p}\e_{0},}
then 
\eq{\label{Main1}d_{\mc{H}}(M,S_{R}(x_{M}))\leq \fr{c^{\a}R}{\|H\|_{p}^{\a}}\|\mathring{A}\|_{p}^{\a}\equiv R\e^{\a}}
 and $M$ is $\e^{\a}$-quasi-isometric to a sphere $S_{R}$ with a certain radius $R.$
}

\Theo{rem}{}{
(i)~By $\e^{\a}$-quasi-isometric we mean that a suitable diffeomorphism $F$ from $M$ into $S_{R}$ satisfies
\eq {\left|d(F(x_{1}),F(x_{2}))-d(x_{1},x_{2})\right|\leqslant R\e^{\a}}
for any $x_{1},x_{2}\in M.$

(ii)~The radius $R$ can be expressed in terms of $\|H\|_{p},$ compare \cite[Cor.~4.6]{RothScheuer:10/2015}.

(iii)~The assumption $|M|=1$ is only for simplification. By scaling it is easy to obtain a scale-invariant version for arbitrary volume.
}

In \cref{ConfFlat} we generalize this theorem to conformally flat ambient spaces. 

The history of the problem to control the closeness to a sphere by curvature quantities is quite long, starting from the well known \it{Nabelpunktsatz}. We refer to the bibliography in \cite{Perez:/2011} for a quite detailed overview. Let us only mention several results which have appeared recently. For surfaces, $n=2,$ a quite straightforward calculation due to Andrews yields an explicit $C^{0}$-estimate for convex hypersurfaces, cf. \cite[Prop.~4, Lemma~5]{Andrews:/1999},
\eq{\label{Andrews}\left|\inpr{x-q}{\nu}-\fr{1}{8\pi}\int_{M}H\right|\leq C|M|\|\mathring{A}\|_{\infty},}
where $x$ is the embedding vector and $q$ is the Steiner point.
In \cref{Optimality} we use the inverse mean curvature flow (IMCF) in the hyperbolic space to prove that the power on the right-hand side of \eqref{Andrews} can not be improved to $\|\mathring{A}\|_{\infty}^{\a},$ $\a>1,$ which is in turn then not possible either in \cref{Main}. The latter proof relies on a recent example due to Hung and Wang, \cite[Thm.~1, Prop.~5]{HungWang:12/2014}, that the convergence after rescaling in the IMCF can not be too fast in the hyperbolic space.

 For strictly convex hypersurfaces of $\R^{n+1}$ there is the following estimate of circumradius $R$ minus inradius $r$ due to Leichtwei\ss, cf. \cite[Thm.~1.4, eq.~(38)]{Leichtweis:08/1999}:
\eq{R-r\leq c_{n}\max\limits_{x\in M}(R_{n}(x)-R_{1}(x)), }
where $R_{1}\leq\dots\leq R_{n}$
are the ordered radii of curvature.
 \cref{Main} deals with estimates in dependence of integral pinching. For the case $n=2,$ an estimate similar to \eqref{Main1} with a better constant was obtained by De Lellis and M{\"u}ller, cf. \cite{De-LellisMuller:/2006}

In \cite[Cor.~1.2]{Perez:/2011} Perez derived a \it{qualitative} solution and obtained under certain assumptions, for given $\e>0,$ a $\delta>0,$ such that
\eq{\|\mathring{A}\|_{p}<\d}
implies
\eq{\label{Intro1}d_{\mc{H}}(M,S_{r_{0}}(x))<\e.}
In \cite[p.~xvi]{Perez:/2011} the author posed the derivation of an explicit $\d$ as a question of interest.

Note that in \eqref{Main1} we did not achieve a constant independent of the size of the curvature itself. The constant is only uniform in the class of hypersurfaces with a fixed bound on the curvature of the hypersurface.

The following theorem, due to Grosjean and the first author, \cite[Thm.~1.4]{GrosjeanRoth:/2012},
already provides this conclusion, however only with the additional assumption of smallness of the oscillation of the mean curvature itself:

\Theo{thm}{Roth}{\cite[Thm.~1.4]{GrosjeanRoth:/2012}\\
Let $(M^{n},g)$ be a compact, connected and oriented $n$-dimensional Riemannian manifold without boundary isometrically immersed by $\phi$ in $\R^{n+1}.$ Let $\e<1,$ $r,q>n,$ $s\geq r$ and $c>0.$ Let us assume that $|M|^{\fr 1n}\|H\|_{q}\leq c.$ Then there exist positive constants $C=C(n,q,c),$ $\a=\a(q,n),$ such that if $\e^{\a}\leq \fr 1C,$
\eq{\label{1}\|\mathring A\|_{r}\leq \|H\|_{r}\e}
and
\eq{\label{2}\|H^{2}-\|H\|_{s}^{2}\|_{\fr r2}\leq \|H\|_{r}^{2}\e,}
then $M$ is $\e^{\a}$-Hausdorff close to $S_{\fr{1}{\|H\|_{2}}}(x_{M}).$ Moreover if $|M|^{\fr 1n}\|A\|_{q}\leq c,$ then $M$ is diffeomorphic and $\e^{\a}$-quasi-isometric to $S_{\fr{1}{\|H\|_{2}}}(x_{M}).$
}

Note that in this theorem, $L^{p}$-norms are defined slightly different, namely such that the $L^{p}$-norms of scale-invariant functions are scale-invariant. Our notation corresponds to the one in \cite{Perez:/2011}. This ambiguity does not cause any problems, since we prove \cref{Main} for $|M|=1.$ Also note the typo in \cite[Thm.~1.4]{GrosjeanRoth:/2012}, where the $\a$ is missing in the conclusion.

In \cite[Thm.~3.1]{Roth:/2015}, which also covers other ambient spaces, \eqref{2} was replaced by an assumption on the gradient of $H.$ However, with the help of the following theorem due to Perez it is possible to get rid of \eqref{2} completely.

\Theo{thm}{Perez}{\cite[Thm.~1.1]{Perez:/2011}\\
Let $p>n\geq 2$ and $c_{0}>0$ be given. Then there is a constant $C>0,$ depending only on $n,$ $p$ and $c_{0},$ such that:\\
If $\Sigma\sub\R^{n+1}$ is a smooth, closed and connected $n$-dimensional hypersurface with 
\eq{|\Sigma|=1}
and 
\eq{\|A\|_{p}\leq c_{0},}
then 
\eq{\min_{\l\in\R}\|A-\l g\|_{p}\leq C\|\mathring{A}\|_{p}.}
}

The proof of \cref{Main} is a combination of \cref{Roth} and \cref{Perez}.

\section{Proofs of \cref{Main}}

{\bf{Proof no.~1:}}
Without loss of generality we may suppose that $M$ is of class $C^{\infty},$ since both sides of the inequality are continuous with respect to the $C^{2}$-norm and hence the general result can then be achieved by approximation.

Using \cref{Perez}, we obtain a $\l_{0}\in\R,$ such that
\eq{\label{Proof1.1}\|A-\l_{0}g\|_{p}\leq C'\|\mathring{A}\|_{p},}
where $C'=C'(n,p,\|A\|_{p}).$
Let us calculate
\eq{\label{Proof1.2}\|H^{2}-\|H\|_{p}^{2}\|_{\fr p2}&\leq \|H^{2}-\l_{0}^{2}\|_{\fr p2}+\|\l_{0}^{2}-\|H\|^{2}_{p}\|_{\fr p2}\\
						&=\(\int_{M}|H-\l_{0}|^{\fr p2}|H+\l_{0}|^{\fr p2}\)^{\fr 2p}+|\l_{0}^{2}-\|H\|^{2}_{p}|\\
						&\leq 2(\|H\|_{p}+|\l_{0}|)\|H-\l_{0}\|_{p}\\
						&\leq c_{n}(\|H\|_{p}+|\l_{0}|)\|A-\l_{0}g\|_{p}\\
						&\leq c'\|H\|_{p}\|\mathring{A}\|_{p},}
where $c'=c'(n,p,\|A\|_{p}).$ The last inequality is due to the fact that
\eq{|\l_{0}-\|H\|_{p}|\leq c''\|\mathring{A}\|_{p}.}

Defining
\eq{c=\max(1,c'),}
\eq{\e=\fr{c\|\mathring{A}\|_{p}}{\|H\|_{p}},}
and
\eq{\e_{0}:=\fr{\min\(1,C^{-\fr{1}{\a}}\)}{2c}} then by \eqref{Main2},
\eq{\e\leq c\e_{0}=\fr 12\min\(1,C^{-\fr{1}{\a}}\), }
where $\a$ and $C$ are the constants from \cref{Roth}. Furthermore we have
\eq{\|\mathring A\|_{p}\leq \|H\|_{p}\e}
 and
 \eq{\|H^{2}-\|H\|_{p}^{2}\|_{\fr p2}\leq \|H\|_{p}^{2}\e.}
 Thus we may apply \cref{Roth} to conclude that $M$ is $\e^{\a}$-close to a sphere. 
 
 The proof of the theorem we applied here, \cref{Roth}, relies on a pinching result for the first eigenvalue which was proven in \cite{GrosjeanRoth:/2012} for a much more general class of ambient spaces. Thus it might not be easily accessible from our point of view. For convenience we want to repeat their main steps of the proof of this theorem in our Euclidean setting, see \cite[p.~487]{GrosjeanRoth:/2012} for the original one. For this purpose we use a recent pinching result for the first eigenvalue of the Laplace operator by both of the authors, cf.~\cite[Thm.~1.1]{RothScheuer:10/2015}. This, and also the original proof in \cite{GrosjeanRoth:/2012}, uses the fact that pinching of the Ricci tensor can be controlled by pinching of the traceless second fundamental form. Then we apply an eigenvalue pinching result due to Aubry, which was proved in \cite[Prop.~1.5]{Aubry:/2007} and can also be found in \cite[Thm.~1.6]{Aubry:/2009}. It says that for $p>n/2$,
a complete Riemannian manifold $(M^{n},g)$ with
\eq{
\frac{1}{|M|}\int_{M}(\underline{\mathrm{Ric}}-(n-1))_{-}^{p}<\frac{1}{C(p,n)}
}
is compact and satisfies
\eq{
\lambda_{1}\geqslant n\left(1-C(n,p)\left(\frac{1}{|M|}\int_{M}(\underline{\mathrm{Ric}}-(n-1))_{-}^{p}\right)^{\frac 1p}\right),
}
where $\underline{\mathrm{Ric}}$ denotes the smallest eigenvalue of the Ricci tensor and $x_{-}=\max(0,-x).$

{\bf{Proof no.~2:}}
Due to the Gauss equation and a simple calculation we obtain a formula for the Ricci tensor in terms of the second fundamental form, namely we obtain
\eq{R_{ij}-(n-1)H^{2}g_{ij}&=(n-2)H(h_{ij}-Hg_{ij})-(h_{ik}-Hg_{ik})(h^{k}_{j}-H\d^{k}_{j}).}
Thus
\eq{\|\mrm{Ric}-(n-1)\|H\|^{2}_{p}g\|_{\fr p2}&\leq c\|H\|_{p}\|\mathring{A}\|_{p}+c\|\mathring{A}\|_{\fr p2}^{2}+\|H^{2}-\|H\|_{p}^{2}\|_{\fr p2}\\
					&\leq c\|H\|_{p}\|\mathring{A}\|_{p},}
where we used \eqref{Proof1.2} and $c=c(n,p,\|A\|_{p}).$
 Using a scaled version of Aubry's eigenvalue estimate we obtain the existence of a constant $\e_{0}=\e_{0}(n,p,\|A\|_{p}),$ such that 
\eq{\|\mathring{A}\|_{p}\leq \e_{0}\|H\|_{p}}
implies
\eq{\l_{1}&\geq n\(\|H\|_{p}^{2}-c\|\mrm{Ric}-(n-1)\|H\|_{p}^{2}g\|_{\fr p2}\)\\
		&\geq n\|H\|_{p}^{2}-c\|H\|_{p}^{2}\fr{\|\mathring{A}\|_{p}}{\|H\|_{p}}\\
		&\geq \(1-c\fr{\|\mathring{A}\|_{p}}{\|H\|_{p}}\)n\|H\|_{p}^{2}.}
Now we can apply the abstract eigenvalue pinching result \cite[Thm.~1.1]{RothScheuer:10/2015}, applied to the tensors $S=T=\id.$

\section{Generalization to conformally flat manifolds}\label{ConfFlat}

Using that the property of a hypersurface to be totally umbilic is invariant with respect to a conformal change of the ambient metric, we easily obtain the following generalization to conformally flat manifolds, which in particular include the half-sphere and the hyperbolic space and improves the $\e^{\a}$-proximity statement in \cite[Thm.~1.3]{GrosjeanRoth:/2012} in the sense that it removes an assumption similar to \eqref{2}.

\Theo{thm}{Cor}{Let $\O\sub\R^{n+1}$ be open and let $N^{n+1}=(\O,\-g)$ be a conformally flat Riemannian manifold, i.e.
\eq{\label{Conf}\-g=e^{2\psi}\~g,}
where $\~g$ is the Euclidean metric and $\psi\in C^{\infty}(\O).$ Let $M^{n}\hra N^{n+1}$ be a closed, connected, oriented and immersed $C^{2}$-hypersurface. Let $p>n\geq 2.$
Then there exist constants $c$ and $\e_{0},$ depending on $n,$ $p,$ $|M|,$ $\|\~A\|_{p}$ and $\|\psi\|_{\infty,M},$ as well as a constant $\a=\a(n,p),$ such that whenever there holds
\eq{\|\mathring{A}\|_{p}\leq \|\~H\|_{p}\e_{0},}
there also holds
\eq{d_{\mc{H}}(M,S_{R})\leq \fr{cR}{\|\~H\|_{p}^{\a}}\|\mathring{A}\|_{p}^{\a},}
 where $S_{R}$ is the image of a Euclidean sphere considered as a hypersurface in $N^{n+1},$ $\|\~A\|_{p}$ and $\|\~H\|_{p}$ are the corresponding Euclidean quantities and the Hausdorff distance is measured with respect to the metric $\-g.$
}

\Theo{rem}{}{
Since in conformally flat spaces the scaling behaviour of the second fundamental form with respect to homotheties heavily depends on the nature of the ambient space, in this case there seems to be no way to give a general scale invariant estimate. This is the reason why this closeness estimate is only uniformly valid in the class of $C^{2}$-bounded hypersurfaces.

Furthermore note that for example in all simply connected space forms the hypersurface $S_{R}$ is actually a geodesic sphere. This follows from the fact that in those spaces totally umbilical hypersurfaces are spheres and total umbilicity invariant with respect to conformal transformations of the ambient space, as will be apparant from the following proof of \cref{Cor}. 

Thus \cref{Cor} gives an explicit spherical closeness estimate of almost-umbilical hypersurfaces in the hyperbolic space as well as in the half-sphere of constant positive sectional curvature. 
}

\pf{
Under a conformal relation of the metrics as in \eqref{Conf} the corresponding induced geometric quantities of the the embedded hypersurface $M$ are related as follows.
\eq{\label{Conf1}g_{ij}=e^{2\psi}\~g_{ij}}
and
\eq{\label{Conf2}h_{ij}e^{-\psi}=\~h_{ij}+\psi_{\b}\~\nu^{\b}\~g_{ij},}
where $\~\nu$ is the normal to $M.$
Those formulae can be found in \cite[Prop.~1.1.11]{Gerhardt:/2006}.
Hence
\eq{\label{Conf3}h_{ij}-Hg_{ij}=e^{\psi}(\~h_{ij}-\~H \~g_{ij})}
and hence
\eq{c\|\mathring{\~A}\|_{p}\leq\|\mathring A\|_{p}\leq C\|\mathring{\~A}\|_{p},}
where the constants depend on $\|\psi\|_{\infty,M}.$
Since the Euclidean and the conformal Hausdorff distances are equivalent whenever $|\psi|$ is bounded, we obtain the result after applying \cref{Main}.
}

Due to a well known interpolation theorem for convex hypersurfaces of Riemannian manifolds we obtain the following gradient stability estimate in space forms.

\Theo{cor}{}{
Let $N^{n+1}$ be the Euclidean space, the hyperbolic space or the sphere. Let $M$ as in \cref{Cor} be additionally strictly convex, where we also assume that $\-g$ is given in geodesic polar coordinates
\eq{\-g=dr^{2}+\vt^{2}(r)\s_{ij}dx^{i}dx^{j}\equiv dr^{2}+\-g_{ij}dx^{i}dx^{j}}
with suitable $\vt$ depending on the space form. Let $p>n.$ Then there exist constants $c$ and $\e_{0}$ depending on $n,$ $p,$ $|M|,$ $\|\~A\|_{p}$ and $\|\psi\|_{\infty},$ as well as a constant $\a=\a(p,n),$
such that
\eq{\label{Grad1}\|\mathring{A}\|_{p}\leq \|\~H\|_{p}\e_{0}}
implies
\eq{v=\sqrt{1+\-g^{ij}u_{i}u_{j}}\leq e^{\fr{cR}{\|\~H\|_{p}^{\a}}\|\mathring{A}\|_{p}^{\a}},}
where 
\eq{M=\{(x^{0},x^{i})\cn x^{0}=u(x^{i}), (x^{i})\in\mc{S}_{0}\}}
is a suitable graph representation over a geodesic sphere $\mc{S}_{0}\hra N^{n+1}$ and $(\-g^{ij})$ is the inverse of $(\-g_{ij}).$
}

\pf{It is well known that a strictly convex hypersurface of $\S^{n+1}$ is contained in an open hemisphere, cf. \cite{CarmoWarner:/1970} for the smooth case and also \cite[Cor.~1.2]{MakowskiScheuer:/2013} for the $C^{2}$-case. Thus $M$ is covered by a conformally flat coordinate system as in \cref{Cor}, which is thus applicable. Let $\mc{S}_{0}$ be the corresponding sphere with center $x_{M},$ then we can write $M$ as a graph over $\mc{S}_{0}$ due to the strict convexity.
Thus we may apply the well-known interpolation estimate
\eq{v\leq e^{\-\k\osc u},}
cf. \cite[Thm.~2.7.10]{Gerhardt:/2006}, where 
\eq{\osc u=\max u-\min u}
and where $\-\k$ is a lower bound for the principal curvatures of the coordinate slices $\{r=\mrm{const}\}.$ The latter, however, only depends on $\|\psi\|_{\infty}$ as well.
}

\section{An optimality result}\label{Optimality}

We prove the optimality of the estimate \eqref{Andrews} in the sense that there is no hope to derive a uniform estimate of the form
\eq{\label{FakeMain1}d_{\mc{H}}(M,S_{R}(x_{0}))\leq c\|\mathring{A}\|_{\infty}^{\a},\quad \a>1,}
in the class of uniformly $C^{\infty}$-bounded hypersurfaces $M.$
To be precise, for $\a>1$ we get the following negation of \eqref{FakeMain1} in the class of uniformly convex hypersurfaces and for all $n\geq 2.$

\Theo{thm}{Opt}{
Let $n\geq 2$ and $C=2\max(|S_{2}(0)|,\|\-A_{S_{2}}\|_{\infty}).$ For all $\a>1$ and for all $k\in\N$ there exists a uniformly convex smooth hypersurface $M_{k}\hra \R^{n+1}$ with 
\eq{\max(\|A_{k}\|_{\infty},|M_{k}|)\leq C,} such that 
\eq{\label{Opta}\|\mathring{A}_{k}\|_{\infty}<\fr 1k}
and for all spheres $S\sub\R^{n+1}$ there holds
\eq{\label{Optb}\quad d_{\mc{H}}(M_{k},S)>k\|\mathring{A}_{k}\|_{\infty}^{\a}.}
Here $\-A_{S_{2}}$ denotes the second fundamental form of the sphere with radius $2.$}

 In a recent paper, Drach gave a counterexample to an improved spherical closeness estimate in the class of $C^{1,1}$ hypersurfaces, namely a special spindle shaped hypersurface, cf. the construction at the beginning of \cite[Sec.~2]{Drach:/2015} and also compare cf.~\cite[Thm.~1]{Drach:/2015}. However, since we consider \eqref{Main1} in the space of at least $C^{2}$-hypersurfaces, we need to find a different contradiction to \eqref{FakeMain1}. This contradiction is deduced along the inverse mean curvature flow in the hyperbolic space.

Before we prove \cref{Opt}, let us for convenience recall the relevant facts about the inverse mean curvature flow in the hyperbolic space $\H^{n+1}.$
There one considers a time parameter family of embeddings of closed, starshaped and mean-convex hypersurfaces
\eq{x\cn [0,T^{*})\times M\hra \H^{n+1},}
which solves
\eq{\dot{x}=\fr{1}{H}\nu,}
where $H=g^{ij}h_{ij}$ and $\nu$ is the outward unit normal to $M_{t}=x(t,M).$ Note that we have switched the notation of $H$ in this context due to a better comparability with the literature. It is known, cf. \cite[Lemma~3.2]{Gerhardt:11/2011}, that for an initial starshaped and mean-convex hypersurface $M_{0}$ the flow exists for all times and all the flow hypersurfaces can be written as a graph over a fixed geodesic sphere $\mc{S}_{0},$ 
\eq{M_{t}=\{(x^{0},x^{i})\cn x^{0}(t,\xi)=u(t,x^{i}(t,\xi))\},}
where $u$ describes the radial distance to the center of $\mc{S}_{0}.$ In \cite[Thm.~1.2]{Gerhardt:11/2011} Gerhardt claimed to have shown convergence of the rescaled hypersurfaces 
\eq{\hat M_{t}=\mrm{graph}~\hat u\equiv\mrm{graph}\(u-\fr{t}{n}\)}
to a geodesic sphere. However, as was pointed out in \cite[Thm.~1]{HungWang:12/2014} with the help of a concrete counterexample, the limit function of $\hat u$ is not constant in general. In particular the authors proved that there is a starshaped and mean-convex initial hypersurface $M_{0},$ such that the limit hypersurface is not of constant curvature, in particular not a geodesic sphere.
However, there is a smooth limit function to which the $\hat M_{t}$ converge smoothly, compare the proof of \cite[Thm.~6.11]{Gerhardt:11/2011} and also compare \cite[Thm.~1.2]{Scheuer:05/2015}.

In order to relate the convergence results of the IMCF in the hyperbolic space with the rigidity estimate \eqref{Main1} in the Euclidean space, we have to look at the hyperbolic flow in the conformally flat model. In \cite{Gerhardt:11/2011} the Poincar\'e ball model in the ball of radius $2$ was considered.
Let $r$ denote the geodesic distance to the center of $\mc{S}_{0}$ in $\H^{n+1},$ then the by the coordinate change
\eq{\rho=2-\fr{4}{e^{r}+1}}
the representation of the hyperbolic metric transforms like
\eq{\-g=dr^{2}+\sinh^{2}(r)\s_{ij}dx^{i}dx^{j}=\fr{1}{\(1-\fr{1}{4}\rho^{2}\)^{2}}(d\rho^{2}+\rho^{2}\s_{ij}dx^{i}dx^{j})\equiv e^{2\psi}\~g,}
where $\s_{ij}$ is the standard round metric of the sphere $\mc{S}_{0}.$
Then the convergence 
\eq{u-\fr{t}{n}\ra \hat u_{\infty}}
in the original coordinates is equivalent to the convergence of
\eq{\label{Opt1}(2-w)e^{\fr tn}\ra \hat w_{\infty},}
where 
\eq{\label{Opt2}w=2-\fr{4}{e^{u}+1}}
 and where $\hat w_{\infty}$ is a strictly positive function due to \cite[Lemma~3.1]{Gerhardt:11/2011}.  

The proof of \cref{Opt} is very similar to the proof of a corresponding positive result in this direction by the second author. In \cite{Scheuer:07/2015} he proved that due to a strong decay of the traceless second fundamental form along the IMCF in $\R^{n+1}$ we indeed obtain spherical roundness in this case without rescaling. The idea how to obtain a negative result in the hyperbolic space is that if we could improve the spherical closeness, then we could mimic the proof in \cite{Scheuer:07/2015} to deduce a roundness result in $\H^{n+1},$ which is not possible in view of Hung's and Wang's paper.
 
The idea of the proof of \cref{Opt} goes as follows: The estimate in \eqref{Opt1} provides closeness of the flow hypersurfaces to the sphere of radius $2$ in the ball model. The order of the closeness is $e^{-\fr tn}.$ The traceless second fundamental form decays correspondingly, as we will point in more detail later in the proof. But if we had this additional exponent $\a$ in the spherical closeness estimate, we could even deduce better spherical closeness (to a sphere different from $S_{2}$) than we have in \eqref{Opt1} and then we would be able to translate this to a spherical closeness in the hyperbolic space. This would in turn yield a contradiction to Hung's and Wang's result. Now let us prove \cref{Opt} in detail. First we need some helpful notation and an auxilliary result. 

\Theo{defn}{Starshaped}{
(i)~Let $N$ be either the Euclidean space, the hyperbolic space or an open hemisphere. For a starshaped hypersurface $M\hra N,$ let
$M^{*}$ be the set of points in $N,$ with respect to which $M$ is starshaped.

(ii)~For a starshaped hypersurface $M\hra N$ let $p\in M^{*}.$ Then for the graph representation
\eq{M=\{(r,x^{i})\cn r=u(x^{i}), (x^{i})\in \mc{S}_{p}\},}
by 
\eq{\osc_{p}u=\max\limits_{x\in \mc{S}_{p}}u(x)-\min\limits_{x\in\mc{S}_{p}}u(x)}
we denote the oscillation of the geodesic distance of the point $(u,x^{i})$ to the point $p.$
Here $\mc{S}_{p}$ denotes a geodesic sphere around $p.$  
}

By a simple argument we obtain the following alternative for a general expanding sequence of hypersurfaces with controlled oscillation.

\Theo{lemma}{Alternative}{
Let $N$ be as in \cref{Starshaped} and $M_{t}\hra N,$ $0\leq t\in\R,$ be a family of starshaped hypersurfaces such that 
\eq{M^{*}_{t}\sub M^{*}_{s}\quad\forall s\geq t} 
and such that for each $\tau_{0}\geq 0$ and $p\in M_{\tau_{0}}^{*}$ there exists a constant $c,$ such that for all $t_{0}\geq \tau_{0},$
\eq{\osc_{p}u_{t}\leq c\osc_{p}u_{t_{0}}\quad\forall t\geq t_{0}.} 
Then for fixed $p,$ $\osc_{p}u_{t}$ does not have zero as a limit value for $t\ra\infty$ unless
\eq{\osc_{p}u_{t}\ra0,\quad t\ra\infty.}
}

\pf{For given $\e>0,$ if zero is a limit point, we may choose $t_{0},$ such that
\eq{\osc_{p}u_{t_{0}}\leq \fr{\e}{c},}
then
\eq{\osc_{p}u_{t}\leq c\osc_{p}u_{t_{0}}\leq \e\quad\forall t\geq t_{0}.}
}

Now we can prove \cref{Opt}.

\pf{
Assume the contrary, i.e. that there exists $\a>1$ and $k\in\N,$ such that for all uniformly convex hypersurfaces $\~M\hra\R^{n+1}$ with
\eq{\max(|\~M|,\|\~A\|_{\infty})\leq C}
we have that
\eq{\|\mathring{\~A}\|_{\infty}<\fr{1}{k}}
implies
\eq{\label{Opt3}\~d_{\mc{H}}(\~M,\~S)\leq k\|\mathring{\~A}\|_{\infty}^{\a}}  
for some suitable sphere $\~S\sub\R^{n+1},$ where the Hausdorff distance is measured with respect to the Euclidean metric.
 According to \cite[Thm.~1]{HungWang:12/2014} for $n=2$ and \cite[Sec.~4]{HungWang:12/2014} for $n\geq 3$ there exists a starshaped and mean-convex hypersurface $M_{0}\hra \H^{n+1},$ such that for no graph representation 
\eq{M_{t}=\mrm{graph}~u}
 the rescaled IMCF flow hypersurfaces 
\eq{\hat M_{t}=\mrm{graph}\(u-\fr tn\)\equiv\mrm{graph}~\hat u}
converge to a geodesic sphere.
However, for each graph representation, we obtain smooth convergence of 
\eq{\hat u\ra \hat u_{\infty}. } In \cite[Thm.~1.2~(2)]{Scheuer:05/2015} it is deduced that
\eq{\|\mathring{A}\|_{\infty}\leq ce^{-\fr{2t}{n}},}
where $c=c(n,M_{0}).$ 
Now fix a graph representation around $p\in M_{0}^{*}.$
From \eqref{Conf3} we obtain that the corresponding Euclidean traceless part decays like
\eq{\|\mathring{\~A}\|_{\infty}=\|e^{\psi}\mathring{A}\|_{\infty}\leq e^{\psi}_{\max}e^{-\fr{2t}{n}}, }
where
\eq{e^{\psi}_{\max}=\fr{1}{\(1-\fr{1}{4}w_{\max}^{2}\)}}
with $w$ as in \eqref{Opt2} and
\eq{w_{\max}=\max\limits_{x\in \mc{S}_{p}}w(x).}
 Due to \eqref{Opt1} we obtain
\eq{\|\mathring{\~A}\|_{\infty}\leq ce^{-\fr tn}}
and due to the $C^{\infty}$-convergence of $w\ra 2,$ we are in the situation to apply our assumption and obtain \eqref{Opt3}, whenever $t$ is large enough. We obtain a sequence of spheres $\~S_{\~R_{t}}\sub \R^{n+1},$
such that
\eq{\~d_{\mc{H}}(\~M_{t},\~S_{\~R_{t}})\leq ce^{-\fr{\a}{n}t}.}
Due to \eqref{Opt1} we even have 
\eq{\~S_{\~R_{t}}\sub B_{2}(0),}
for large times $t.$

Now let us switch back to the hyperbolic space. The spheres $\~S_{\~R_{t}}$ are geodesic spheres in $\H^{n+1}$ as well since total umbilicity is preserved under a conformal transformation and in the Euclidean space as well as in the hyperbolic space for closed and embedded hypersurfaces total umbilicity is tantamount to being a geodesic sphere. We denote these spheres in $\H^{n+1}$ by $S_{R_{t}}.$
For the corresponding hyperbolic Hausdorff distance we deduce
\eq{\label{Opt5}d_{\mc{H}}(M_{t},S_{R_{t}})\leq e^{\psi}_{\max}\~d_{\mc{H}}(\~M_{t},\~S_{\~R_{t}})\leq ce^{\fr{1-\a}{n}t},}
which converges to $0$ as $t\ra\infty.$

Since the inradius of the $M_{t}$ converges to infinity and for large $t$ the $M_{t}$ are strictly convex, for each $\d>0$ we find $t_{0}>0,$ such that
\eq{\-B_{\d}(p)\sub M_{t_{0}}^{*}\sub M_{t}^{*}\quad\forall t\geq t_{0},}
where the latter inclusion is due to the fact that starshapedness around a given point is preserved. 
According to \cite[Prop.~3.2, Lemma~3.5]{Scheuer:05/2015}, there holds for the oscillation of $u$ that for all $\tau_{0},$ all $q\in M_{\tau_{0}}^{*}$ and all $t_{0}\geq\tau_{0}$ we have
\eq{\label{Opt4}\osc_{q}u(t,\cdot)\leq c\osc_{q}u(t_{0},\cdot)\quad\forall t\geq t_{0},}
where $c$ depends on $n$ and on a lower bound on the minimal distance of $q$ to $M_{\tau_{0}}.$
So in particular, if we choose 
\eq{\d=c\osc_{p}u(0,\cdot),}
we find that the oscillation of each $M_{t}$ is minimized within the set $\-B_{\d}(p)\cn$
\eq{\argmin\displaylimits_{q\in M^{*}_{t}}\osc_{q}u(t,\cdot)\in \-B_{\d}(p)\quad\forall t\geq t_{0},}
because outside $\-B_{\d}(p)$ the oscillation is already larger than it is with respect to $p.$

Due to \eqref{Opt5} we obtain
\eq{\label{Opt6}\osc_{q_{t}}u(t,\cdot)=\min\limits_{q\in\-B_{\d}(p)}\osc_{q}u(t,\cdot)\leq ce^{\fr{1-\a}{n}t}\quad\forall t\geq t_{0}.}
Let $t_{k}$ be a sequence of times with $t_{k}\ra\infty.$ Due to the compactness of $\-B_{\d}(p)$ a subsequence of center points converges,
\eq{q_{t_{k}}\equiv q_{k}\ra q\in\-B_{\d}(p),}
where we did not rename the index of the sequence.
Since
\eq{|\osc_{q_{k}}u(t_{k},\cdot)-\osc_{q}u(t_{k},\cdot)|\leq 2\dist(q_{k},q)\quad \forall k\in\N,}
we obtain in view of \eqref{Opt6},
\eq{\osc_{q}u(t_{k},\cdot)\ra 0,\quad k\ra\infty.}
In view of \eqref{Opt4} and the preservation of starshapedness along IMCF the assumptions of \cref{Alternative} are fulfilled. Applying \cref{Alternative}, we obtain that
\eq{\osc_{q}u(t,\cdot)\ra 0,}
in contradiction to the choice of the initial hypersurface.
}

\Theo{rem}{}{
Note that in turn of the proof we even have shown that for given $\a>1$ and $k\in\N$ as in \cref{Opt}, such a counterexample $M_{k}$ satisfying \eqref{Opta} and \eqref{Optb} must actually occur along the inverse mean curvature flow in the conformally flat version of the IMCF in $\H^{n+1}.$ We only used our contrary assumption within this class of flow hypersurfaces.
}

\section{Concluding remark}

We would like to point out that the techniques in \cref{Optimality} might be useful in other situations. Whenever one would like to estimate the closeness to a sphere in comparison with another geometric quantity, e.g. in comparison with eigenvalue pinching of the Laplacian or also in almost-Schur/almost-CMC type estimates, one could determine how this particular geometric quantity behaves along the IMCF and then determine the best possible roundness estimate using the IMCF in $\H^{n+1}.$ It should often be quite straightforward to derive the best possible decay estimate. 

\subsubsection*{Acknowledgements}
We would like to thank Kostiantyn Drach for several valuable comments.

\bibliographystyle{/Users/J_Mac/Documents/Uni/TexTemplates/hamsplain.bst}
\bibliography{/Users/J_Mac/Documents/Uni/TexTemplates/Bibliography}

\providecommand{\bysame}{\leavevmode\hbox to3em{\hrulefill}\thinspace}
\providecommand{\href}[2]{#2}
\begin{thebibliography}{10}

\bibitem{Andrews:/1999}
Ben Andrews, \emph{Gauss curvature flow: the fate of the rolling stones},
  Invent. Math. \textbf{138} (1999), no.~1, 151--161.

\bibitem{Aubry:/2007}
Erwann Aubry, \emph{Finiteness of $\pi_{1}$ and geometric inequalities in
  almost positive {R}icci curvature}, Ann. Sci. l'Ecole Norm. Sup{\'e}r.
  \textbf{40} (2007), no.~4, 675--695.

\bibitem{Aubry:/2009}
\bysame, \emph{Diameter pinching in almost positive {R}icci curvature},
  Comment. Math. Helv. \textbf{84} (2009), no.~2, 223--233.

\bibitem{De-LellisMuller:/2006}
Camillo De~Lellis and Stefan M{\"u}ller, \emph{A ${C}^0$-estimate for nearly
  umbilical surfaces}, Calc. Var. Partial Differ. Equ. \textbf{26} (2006),
  no.~3, 283--296.

\bibitem{CarmoWarner:/1970}
Manfredo Do~Carmo and Frank Warner, \emph{Rigidity and convexity of
  hypersurfaces in spheres}, J. Differ. Geom. \textbf{4} (1970), no.~2,
  133--144.

\bibitem{Drach:/2015}
Kostiantyn Drach, \emph{Some sharp estimates for convex hypersurfaces of
  pinched normal curvature}, J. Math. Phys. Anal. Geom. \textbf{11} (2015),
  no.~2, 111--122.

\bibitem{Gerhardt:/2006}
Claus Gerhardt, \emph{Curvature problems}, Series in Geometry and Topology,
  vol.~39, International Press of Boston Inc., 2006.

\bibitem{Gerhardt:11/2011}
\bysame, \emph{Inverse curvature flows in hyperbolic space}, J. Differ. Geom.
  \textbf{89} (2011), no.~3, 487--527.

\bibitem{GrosjeanRoth:/2012}
Jean-Fran\c cois Grosjean and Julien Roth, \emph{Eigenvalue pinching and
  application to the stability and the almost umbilicity of hypersurfaces},
  Math. Z. \textbf{271} (2012), no.~1-2, 469--488.

\bibitem{HungWang:12/2014}
Pei-Ken Hung and Mu~Tao Wang, \emph{Inverse mean curvature flows in the
  hyperbolic 3-space revisited}, Calc. Var. Partial Differ. Equ. (2014),
  {10.1007/s00526-014-0780-3}, Online first.

\bibitem{Leichtweis:08/1999}
Kurt Leichtwei\ss, \emph{Nearly umbilical ovaloids in the n-space are close to
  spheres}, Result. Math. \textbf{36} (1999), no.~1-2, 102--109.

\bibitem{MakowskiScheuer:/2013}
Matthias Makowski and Julian Scheuer, \emph{Rigidity results, inverse curvature
  flows and {A}lexandrov-{F}enchel type inequalities in the sphere}, to appear
  in Asian J. Math. (2016).

\bibitem{Perez:/2011}
Daniel Perez, \emph{On nearly umbilical hypersurfaces}, Ph.D. thesis, Zuerich,
  2011.

\bibitem{Roth:/2015}
Julien Roth, \emph{A new result about almost umbilical hypersurfaces of real
  space forms}, Bull. Aust. Math. Soc. \textbf{91} (2015), no.~1, 145--154.

\bibitem{RothScheuer:10/2015}
Julien Roth and Julian Scheuer, \emph{Pinching of the first eigenvalue for
  second order operators on hypersurfaces of the {E}uclidean space}, preprint
  available at \href{http://arxiv.org/abs/1510.03722}{arxiv:1510.03722}, 2015.

\bibitem{Scheuer:05/2015}
Julian Scheuer, \emph{Non-scale-invariant inverse curvature flows in hyperbolic
  space}, Calc. Var. Partial Differ. Equ. \textbf{53} (2015), no.~1-2, 91--123.

\bibitem{Scheuer:07/2015}
\bysame, \emph{Pinching and asymptotical roundness for inverse curvature flows
  in {E}uclidean space}, J. Geom. Anal. (2015), 1--17,
  {10.1007/s12220-015-9627-1}, Online first.

\end{thebibliography}

\end{document}